\documentclass{amsart}

\usepackage{amsmath,amsthm,amsopn,amstext,amscd,amsfonts,amssymb,mathrsfs,mathtools}
\usepackage{dsfont}
\usepackage{comment}
\usepackage{xcolor}
\usepackage{float}
\usepackage[active]{srcltx}
\usepackage{graphicx, mathdots}

  \usepackage{mathtools}

\usepackage[active]{srcltx}

\usepackage{graphicx, epsfig, subfig}
\usepackage{lscape}
\usepackage{rotating}
\usepackage{caption}
\usepackage{multicol}
\usepackage{colortbl}
\usepackage{nicematrix}
\usepackage{hyperref}

\newmuskip\pFqskip
\mathchardef\pFcomma=\mathcode`, 

\makeatletter
\def\BState{\State\hskip-\ALG@thistlm}
\makeatother

\def\downbar#1{
\setbox10=\hbox{$#1$}
            \dimen10=\ht10 \advance\dimen10 by 2.5pt
            \ifdim \dimen10<15pt 
               \advance\dimen10 by -0.5pt
               \dimen11=\dimen10
               \advance\dimen10 by 2.5pt
               \lower \dimen11
            \else \lower \ht10 \fi
            \hbox {\hskip 1.5pt \vrule height \dimen10 depth \dp10}}
\def\upbar#1{
\setbox10=\hbox{$#1$}
            \dimen10=\ht10 \advance\dimen10 by \dp10 \advance\dimen10 by 2.5pt
            \ifdim \dimen10<15pt 
                \advance\dimen10 by 2pt \fi
            \raise 2.5pt \hbox {\hskip -1.5pt \vrule height \dimen10}}

\begin{document}
\title[Remark on the eigenvalues of a tridiagonal matrix in biogeography]{Remark on the eigenvalues of a tridiagonal matrix in biogeography}
\author{K. Castillo}
\address{CMUC, Department of Mathematics, University of Coimbra, 3001-501 Coimbra, Portugal}
\email{ kenier@mat.uc.pt}

\date{\today}
\begin{abstract}
The main result proved in  [The eigenvalues of a tridiagonal matrix in biogeography, Appl. Math. Comput. 218 (2011) 195-201; MR2821464] by B. Igelnik and D. Simon is virtually the Sylvester determinant.
\end{abstract}
\maketitle
In \cite{IS}, motivated by an interesting model in biogeography\footnote{See \cite{CZ} for recent results on the eigenvalues of this kind of matrices in general.}, B. Igelnik and D. Simon, with a rather long proof, show that the $n+1$ eigenvalues of 
\begin{align*}
\mathbf{A}_{n+1}=\begin{pmatrix}
-1 & 1/n & & &\\
n/n &   -1 &  \ddots & & \\
& \ddots & \ddots &  (n-1)/n &\\
& & 2/n & -1 & n/n \\
&&& 1/n & -1
\end{pmatrix}
\end{align*}
are $-2\,(n-k+1)/n$ for $k=1,\dots, n+1$. In MR2821464, the reviewer for MathScinet  points out that this follows from a result by P. A. Clement \cite{C59}. But even more can be said, because an elementary proof is well known. The eigenvalues of $\mathbf{A}_{n+1}$ are, up to an affine change of the variable, the roots of the following polynomial, which became known as Sylvester determinant once J. J. Sylvester published a note in Nouvelles Annales de Math\'ematiques in 1854 \cite{S54}:
\begin{align*}
p_{n+1}(X)=\det \begin{pmatrix}
X & 1 & & &\\
n &   X &  \ddots & & \\
& \ddots & \ddots &  (n-1) &\\
& & 2 & X & n \\
&&& 1 & X
\end{pmatrix}.
\end{align*}
In fact, taking $X=n\,(x+1)$ we see that $\det(x\, \mathbf{I}-\mathbf{A}_{n+1})=n^{-(n+1)}\,p_{n+1}(X)$. 
It is worth highlighting that O. Taussky and J. Todd published, in 1991, an interesting expository and historical paper \cite{TT}  where, using elementary row and column operations, they give a rather simple proof of the Sylvester determinant after having noted that
\begin{align*}
p_{n+1}(X)=(X-n) p_{n}(X+1),
\end{align*}
and so 
\begin{align}\label{det}
p_{n+1}(X)=\prod_{k=1}^{n+1} (X+n-2k+2).
\end{align}
(It is clear from \eqref{det} that the eigenvalues of $\mathbf{A}_{n+1}$ are $-2\,(n-k+1)/n$ for $k=1,\dots, n+1$.) Taussky and Todd, who also presented another proof from 1866 due to F. Mazza, wrote that the above arguments come from ``A treatise on the Theory of Determinants" by T. Muir in the edition revised and enlarged by W. H. Metzler. Although they omitted a specific reference concerning the origin of these arguments, the reader should notice that Theorem $576$ with $a=X$, $b=1$, $c=-1$ is the Sylvester determinant (the second example in $578$ is exactly the Sylvester determinant for $n=6$).  They also commented that this determinant appears as an exercise since the original Muir edition, which was published in 1882. We also noted that a more general Sylvester's type determinant was considered by A. Cayley in 1857, see \cite[p.  429]{M}. The connection of $p_{n+1}(X)$ and related determinants with orthogonal polynomials on uniform and non-uniform lattices was nicely explored by R. Askey \cite{A03} and O. Holtz \cite{O03} in 2003 during the 4th International Congress of the International Society for Analysis, Applications and Computation (ISAAC) held at York University, Toronto. Both Askey and Holtz presented a proof of the Sylvester determinant, without references to the works cited above. The first proof given by Askey, in which he does not use orthogonal polynomials and that he adapted from the solution of ``the Painvin determinant"  in  ``Problems in Higher Algebra" by D. K. Faddeev and I. S Sominskii, is the same proposed by Muir and Metzler, and reproduced by Taussky and Todd. In any case, we emphasize that this solution for the Sylvester determinant  is already known to everyone who followed during the first undergraduate year the exercises proposed in ``Problems in Linear Algebra" by I. V. Proskuryakov. In \cite[Problem 399, p. 61]{P}, Proskuryakov asks for the value of the determinant $p_{n+1}(X)$. Below we quote verbatim the suggestion given by Proskuryakov to solve the problem \cite[p. 326]{P}, that is the same given by Faddeev and Sominskii in his book for the Painvin determinant, and in turn the same proposed by L. Painvin himself in 1858 \cite[p. 433]{M}:
\vspace{2mm}
\begin{changemargin}{0.8cm}{0.8cm}
{\em ``To each row add all the following rows, from each column subtract the preceding column and show that $D_n(x)$ $[$$p_{n}(X)$$]$ is the given determinant, then $D_n(x)=(x+n-1)D_{n-1}(x-1)$ $[$$p_{n}(X)=(X+n-1) p_{n-1}(X-1)$$]$".}
\end{changemargin}
\vspace{2mm}
It is worth pointing out that the above suggestion probably\footnote{The oldest Russian edition consulted was the third, which was published in 1966. The prefaces of the second and third edition do not contain any mentions of changes related to this problem.} appears since the first Russian edition in 1955. 

\section*{Acknowledgements}
The author would like to thank R. \'Alvarez-Nodarse for kindly sending a copy of the third Russian edition of Proskuryakov's book.
This work was supported by the Centre for Mathematics of the University of Coimbra-UIDB/00324/2020, funded by the Portuguese Government through FCT/ MCTES.

 \end{document}